\documentclass[a4paper]{article}
\usepackage{amsmath,amsthm,amssymb,xcolor}
\usepackage{tikz, tkz-graph, tkz-berge}
\usepackage{hyperref}
\usepackage[normalem]{ulem}
\usepackage{enumitem}
\usepackage[T1]{fontenc}
\usepackage[plain]{fullpage}
\newtheorem{theorem}{Theorem}
\newtheorem{corollary}[theorem]{Corollary}
\newtheorem{proposition}[theorem]{Proposition}
\newtheorem{lemma}[theorem]{Lemma}
\newtheorem{claim}{Claim}
\theoremstyle{definition}
\newtheorem{problem}[theorem]{Problem}
\newtheorem{conjecture}[theorem]{Conjecture}
\renewcommand{\epsilon}{\varepsilon}
\renewcommand{\phi}{\varphi}

\newcommand{\Ra}{\Rightarrow}
\newcommand{\ra}{\rightarrow}

\newcommand{\ind}[1]{\langle #1 \rangle}

\newcommand{\e}{\mathrm{e}}

\newenvironment{proofclaim}[1][{\it Proof of claim. \hspace{0.066cm}}]
	{\noindent {}{#1}{}}{ \strut\hfill $\lozenge$\vspace{2ex}}

\DeclareMathOperator{\unvd}{unvd}
\DeclareMathOperator{\dist}{dist}

\newcommand{\bigO}{\mathcal{O}}

\renewcommand{\Pr}{\mathbf{Pr}}
\newcommand{\E}{\mathbf{E}}

\renewcommand{\leq}{\leqslant}

\renewcommand{\geq}{\geqslant}

\title{Blow-ups and extensions of trees in tournaments}
\author{
Pierre Aboulker\footnote{DIENS, École normale supérieure, CNRS, PSL University, Paris, France. Research supported by the ANR project DAGDigDec (JCJC) ANR-21-CE48-0012 and by the group Casino/ENS Chair on Algorithmics and Machine Learning.} 
\and Frédéric Havet\footnote{Universit\'e C\^ote d'Azur, CNRS, Inria, I3S, Sophia-Antipolis, France. Research supported by research grant DIGRAPHS ANR-19-CE48-0013.}
\and William Lochet\footnote{LIRMM, Universit\'e de Montpellier, CNRS, Montpellier, France. Research supported by research grant ELiT ANR-20-CE48-0008-01.}
\and Raul Lopes\footnotemark[1] $^,$\footnotemark[3]
\and Lucas Picasarri-Arrieta\footnote{National Institute of Informatics, Japan. Research supported by Japan Science and Technology Agency (JST) as part of Adopting Sustainable Partnerships for Innovative Research Ecosystem (ASPIRE), Grant Number JPMJAP2302.}
\and Clément Rambaud\footnotemark[2]
}
\date{}

\begin{document}

\maketitle

\begin{abstract} 
    A class of acyclic digraphs $\mathcal{C}$ is linearly unavoidable if there exists a constant $c$ such that every digraph $D\in \mathcal{C}$ is contained in all tournaments of order $c\cdot |V(D)|$.
    The class of all acyclic digraphs is not linearly avoidable, and Fox, He, and Widgerson recently showed that this is not even the case for acyclic digraphs with bounded maximum degree.
    On the positive side, Thomason and Häggkvist proved that the class of oriented trees is linearly unavoidable.  
    In this work, we generalize this result to
    acyclic digraphs obtained from an oriented tree by adding at most $k$ vertices, and
    $k$-blow-ups of oriented trees, for every fixed integer $k$.

    More precisely, we show that if $D$ is obtained from an oriented tree $F$ of order $n$ by adding $k$ universal vertices, then $D$ is contained in every tournament of order $2\cdot 3^{(k+1)(2k+1)} \cdot n$; and 
    if $D$ is obtained from $F$ by replacing each vertex $u$ by an independent set $X_u$ of size $k$ and every arc $uv$ by all possible arcs from $X_u$ to $X_v$, then $D$ is contained in every tournament of order $2^{10+18k}k \cdot n$. 
\end{abstract}

\sloppy

\section{Introduction}

A digraph is {\bf $p$-unavoidable}, for an integer $p$, if it is contained in every tournament of order $p$. 
For a digraph $D$, the {\bf unavoidability} of $D$, denoted by $\unvd(D)$, is the smallest integer $p$ such that $D$ is $p$-unavoidable, if such a $p$ exists. 
It is easy and well-known that a digraph is $p$-unavoidable for some integer $p$ if and only if it is acyclic.
Indeed, every digraph containing a directed cycle is not contained in any transitive tournament, and so is not unavoidable.
Moreover, as observed by Erd\H{o}s and Moser~\cite{ErMo64}, an immediate induction and Proposition~\ref{prop:double} shows that the transitive tournament of order $n$ is $2^{n-1}$-unavoidable and thus every acyclic digraph of order $n$ is $2^{n-1}$-unavoidable.

\begin{proposition}\label{prop:double}
    Let $D$ be an acyclic digraph. For every source or sink $x$ of $D$,
    $$\unvd(D) \leq 2\unvd(D-x).$$
\end{proposition}

This is somehow tight in the sense that any upper bound on the unavoidability of acyclic digraphs as a function of their order must be exponential. 
Indeed, considering a random tournament and using the First Moment Method, Erd\H{o}s and Moser~\cite{ErMo64} established $\unvd(TT_n)\geq 2^{(n-1)/2}$. 
Since this seminal paper, only little progress has been made on the unavoidability of transitive tournaments.
Using the Local Lemma, one can slightly improve on the lower bound and show that $\unvd(TT_n)\geq 2^{(n+1)/2}$ when $n$ is large enough.
The unavoidability of small transitive tournaments have been established.
Clearly $\unvd(TT_1)=1$, $\unvd(TT_2)=2$,  $\unvd(TT_3)=4$ (because of the directed $3$-cycle), and  $\unvd(TT_4)=8$ (because of the Paley tournament on $7$ vertices).
Reid and Parker~\cite{RePa70} showed $\unvd(TT_5)=14$,  $\unvd(TT_6)=28$, and Sanchez-Flores~\cite{San98} proved $\unvd(TT_7)=54$.
This result and Proposition~\ref{prop:double} imply   $\unvd(TT_n)\leq 54 \times 2^{n-7}$ for all $n\geq 7$.

The {\bf average degree} of a nonempty digraph $D$ is the average degree of its underlying graph, that is $2\frac{|A(D)|}{|V(D)|}$.
The {\bf maximum average degree} of a digraph is the maximum average degree over all its nonempty subdigraphs.
Similarly to Erd\H{o}s and Moser~\cite{ErMo64}, one can show that the unavoidability of an acyclic digraph is exponential in its maximum average degree.

\begin{proposition}\label{prop:lower-average}
    Let $D$ be an acyclic digraph with maximum average degree $\alpha$. 
    Then $\displaystyle \unvd(D) > 2^{\alpha/2}$.
\end{proposition}
\begin{proof}
    Consider a random tournament $T$ on $p\leq 2^{\alpha/2}$ vertices.
    Let $v_1, \dots , v_n$ be a labelling of the vertices of a subdigraph $H$ of $D$ with average degree $\alpha$.
    For each ordered $n$-uple $(w_1, \dots , w_n)$ of distinct vertices in $T$, the probability that $T\langle \{w_1, \dots , w_n\}\rangle$ contains a labelled copy of $H$ (\textit{i.e.} such that
    $w_iw_j$ is an arc in $T$ whenever $v_iv_j$ is an arc in $H$) is
    $\left (\frac{1}{2}\right )^{|A(H)|} =  \left (\frac{1}{2}\right )^{\alpha n/2}$.
    Thus the expected number of labelled copies of $H$ is 
    $\frac{p!}{(p-n)!}\left ( \frac{1}{2} \right ) ^{\alpha n/2} < \frac{p^{n}}{2^{\alpha n/2}}  = \left (\frac{p}{2^{\alpha/2}}\right )^{n}\leq 1$.
\end{proof}

Proposition~\ref{prop:lower-average} gives a non-trivial lower bound on the unavoidability of dense acyclic digraphs, that is with maximum average degree larger than logarithmic in their order.
However, smaller bounds on the unavoidability are expected for sparse digraphs.
In particular, one can ask which families ${\cal F}$ of digraphs are {\bf linearly unavoidable}, that is such that 
there exists a constant $C$ such that  $\unvd(D) \leq C \cdot |V(D)|$ for every  digraph $D$ in ${\cal F}$,  and which are {\bf polynomially unavoidable}, that is such that 
there exists a polynomial $P$ such that $\unvd(D) \leq P(|V(D)|)$ for every  digraph $D$ in ${\cal F}$. 

Note that by Proposition~\ref{prop:lower-average}, the digraphs $D$ of a linearly unavoidable family must have maximum average degree at most $2\log(|V(D)|)+ \bigO(1)$. 
For every integer $n$ with $n \geq 2$, Linial, Saks and S\'os~\cite{LSS83} constructed an $n$-unavoidable digraph with $n$ vertices and $n\log n - C\cdot n\log\log n$ arcs for some fixed constant $C$. 
Those digraphs form a linearly unavoidable family with logarithmic maximum average degree. 

Several other classes of acyclic digraphs have been proved to be linearly unavoidable. Much attention has been devoted to {\bf oriented paths}, {\bf oriented trees}, and {\bf oriented forests}  which are orientations of paths, trees, and forests respectively.
It started with R\'edei's theorem~\cite{Rede34}
which states that the unavoidability of $\vec{P}_n$, the directed path on $n$ vertices, is exactly $n$.

\begin{theorem}[R\'edei~\cite{Rede34}]\label{thm:redei}
    Every tournament has a Hamiltonian directed path.
\end{theorem}

In 1971, Rosenfeld~\cite{rosenfeldJCTB12} conjectured  that  there is an integer $N>7$ such that $\unvd(P) = |P|$ for every oriented path of order at least $N$. 
Several papers gave partial answers to this conjecture \cite{grunbaumJCTB11,alspachDM34,forcadeDM6,straightCN29} until Rosenfeld's conjecture was verified by Thomason, who proved in~\cite{thomasonTAMS296} Rosenfeld's conjecture for $N = 2^{128}$. Finally, Havet and Thomass\'e~\cite{havetJCTB78} showed that $\unvd(P)=|V(P)|$ for every oriented path $P$ except the antidirected paths of order $3$, $5$, and $7$.

Regarding oriented trees, Sumner (see~\cite{reidSSMH18}) made the following celebrated conjecture.
\begin{conjecture}[Sumner]\label{conj:Sumner}
    Every oriented tree of order $n>1$ is $(2n-2)$-unavoidable.
\end{conjecture}
Since every forest is the subdigraph of a tree, this conjecture is equivalent to the analogue for oriented forests.
If true, Sumner's conjecture would be tight. Indeed, the {\bf out-star} $S^+_n$, which is the digraph on $n$ vertices consisting of a vertex dominating the $n-1$ others, is not contained in the regular tournaments of order $2n-3$. 
The first linear upper bound was given by 
H\"aggkvist and Thomason~\cite{haggkvistComb11}. Following improvements of Havet~\cite{Hav02}, Havet and  Thomass\'e~\cite{havetJGT35}, and El Sahili~\cite{elsahiliJCTB92},  Dross and Havet~\cite{DrHa21} used the notion of median order to prove that every oriented tree of order $n \geq 2$ is $\left \lceil \frac{21}{8} n - \frac{47}{16}  \right \rceil$-unavoidable.
K\"uhn, Mycroft and Osthus~\cite{kuhnJCTB101} proved that
Sumner's conjecture is true for all sufficiently large $n$.
\begin{theorem}[K\"uhn, Mycroft, and Osthus~\cite{kuhnJCTB101}]\label{thm:kuhn_mycroft_osthus}
    There is an integer $n_0$ such that for every $n \geq n_0$, every oriented tree of order $n$ is $(2n-2)$-unavoidable.
\end{theorem}

Motivated by those results and a work of Lee~\cite{Lee2017}, Buci\'c, Letzter and Sudakov~\cite{BuLeSu19}
asked whether the digraphs with bounded maximum average degree are linearly unavoidable.
Dragani\'c \textit{et al.}~\cite{draganicCPC30} proved that this is the case  for $k$-th powers of directed paths.
The {\bf $k$-th power} of a digraph $D$ is the digraph, denoted by $D^k$, with vertex set $V(D)$ in which $uv$ in an arc if and only if $\dist_D(u,v)\leq k$.
Fox, He, and Widgerson~\cite{Fox2024} showed that this is not the case in general:
they proved that for any $\Delta \geq 2$ and any sufficiently large $n$, there exists an acyclic digraph $D$ with $n$ vertices and maximum degree $\Delta$ such that $$\unvd(D)\geq n^{\Omega(\Delta^{2/3}/\log^{5/3}\Delta )}.$$
This result raises the two following questions.
Firstly, are the digraphs with bounded maximum average degree polynomially unavoidable? 

\begin{problem}
    Let $\alpha$ be a positive real number.
    Does there exist a polynomial $P_{\alpha}$
    such that $\unvd(D) \leq P_{\alpha}(|V(D)|)$ for every digraph with maximum average degree at most $\alpha$?
\end{problem} 

Secondly, which families of acyclic digraphs (with bounded maximum average degree) are linearly unavoidable?
In this direction,
Fox, He, and Widgerson~\cite{Fox2024} showed that for every acyclic digraph $D$ of maximum degree $\Delta$, if there is a partition $S_1, \dots, S_h$ of $V(D)$
such that for every arc $uv$ of $D$, there exists $i\in [h-1]$ such that $u \in S_i$ and $v \in S_{i+1}$,
then $\unvd(D) \leq h^{10 \Delta \log \Delta} |V(D)|$.

\subsection*{Our results}
In this paper, we give some more partial answers to the above second question by proving some new families of acyclic digraphs to be linearly unavoidable. We are in particular interested in operations to form linearly unavoidable families from others.

\medskip

Let $k$ be a positive integer.
The {\bf $k$-blow-up} of a digraph $D$, denoted by $D[k]$, is obtained by blowing-up each vertex into an independent set of size $k$. Formally, $D[k]$ is defined by 
\begin{eqnarray*}
    V(D[k]) & = & \{ (v,i) \mid v\in V(D), i\in [k]\}, \mbox{ and }\\
    A(D[k]) & = & \{ (u,i)(v,j) \mid uv\in A(D), i\in [k], j\in [k], i \neq j\} .
\end{eqnarray*}
We believe that if a family $\cal F$ is linearly unavoidable, then the family ${\cal F}[k] = \{D[k] \mid D \in {\cal F}\}$ of its $k$-blow-ups is also linearly unavoidable.

\begin{conjecture}\label{conj:blow-up}
    If $\cal F$ is linearly unavoidable and has bounded maximum average degree, then ${\cal F}[k]$ is also linearly unavoidable.
\end{conjecture}

Observe first that the bounded maximum average degree condition in the above conjecture is necessary. 
Indeed, for every $n$, let $D_n$ be the digraph of order $n$ which is the disjoint union of a transitive tournament of order $\lfloor \log n \rfloor + 1$ 
and $n - \lfloor \log n \rfloor -1$ isolated vertices. Each $D_n$ is clearly $n$-unavoidable, but
$D_n[k]$ has maximum average degree $k\lfloor \log n \rfloor$ and so
$\unvd(D_n) \geq 2^{k\lfloor \log n \rfloor /2}\geq \frac{1}{2}n^{k/2}$ by Proposition~\ref{prop:lower-average}.

Observe that the $k$-blow-up of the directed path of order $n$ is contained in the $2k$-th power of the directed path of order $kn$. 
Thus, the result of Dragani\'c \textit{et al.}~\cite{draganicCPC30} implies Conjecture~\ref{conj:blow-up} for directed paths. 
As a more general evidence to this conjecture, we show in Section~\ref{sec:tree-blow-up} that it holds for oriented trees (and thus more generally oriented forests). 
Precisely, we prove in Theorem~\ref{thm:blowup} that the $k$-blow-up of an oriented tree of order $n$ is $(2^{10+18k} \cdot k n)$-unavoidable.
This upper bound is certainly not tight. However, the maximum average degree of the $k$-blow-up of a tree of order $n$ is $2k - 2k/n$. 
Therefore, by Proposition~\ref{prop:lower-average}, the optimal bound is of the form $2^{\Theta(k)} \cdot kn$. Hence we are left with the following question.

\begin{problem}\label{pb:blowup}
    What is the infimum of all the constants $C$ such that for every large enough integer $n$, 
    the $k$-blow-up of every oriented tree of order $n$ is $(C^k \cdot kn)$-unavoidable?
\end{problem}

Generalizing Proposition~\ref{prop:double}, we believe that adding a vertex cannot affect too much the unavoidability, in a sense that it can at most  multiply it by a constant factor.

\begin{conjecture}\label{conj:foisC}
There exists a constant $C$ such that $\unvd(D) \leq C\cdot \unvd(D-v)$ for every acyclic digraph $D$ and for every vertex $v$ of $D$.
\end{conjecture}

Let $A$ be an acyclic digraph and let $k$ be a positive integer.
An acyclic digraph $D$ is a {\bf $k$-extension} of $A$ if there exists a set $S$ of $k$ vertices such that $D-S=A$.
The set $S$ is the {\bf added set} and its vertices are the {\bf added vertices}.
A vertex is {\bf universal} in a digraph $D$ if $N^-(v)\cup N^+(v)\cup\{v\}=V(D)$.
A $k$-extension $D$ of $A$ is {\bf full}, if its added vertices are universal in $D$.
A $k$-extension $D$ of $A$ is a {\bf $k$-twin extension}, if all added vertices have the same in- and out-neighbourhood in $V(D-S)$.
Conjecture~\ref{conj:foisC}, together with Conjecture~\ref{conj:Sumner}, implies the following.

\begin{conjecture}\label{conj:k-ext-constant}
 There is an absolute constant $C$ such that for every integer $k$, every $k$-extension of an oriented tree of order $n$ is $C^{k}(2n-2)$-unavoidable.  
\end{conjecture}

More generally, Conjecture~\ref{conj:foisC} would imply the following.

\begin{conjecture}\label{conj:k-ext}
 If $\cal F$ is linearly unavoidable, then the family of $k$-extensions of digraphs in $\cal F$ is also linearly unavoidable.
\end{conjecture}

In Section~\ref{sec:extension-tree}, we prove this last conjecture for the family of oriented trees (and thus also of oriented forest). More precisely; we show in Corollary~\ref{cor:k-extension-of-trees} that every $k$-extension of a sufficiently large tree $F$ is $\left ( 2 \cdot 3^{\binom{2k+2}{2}} \cdot |V(F)|\right )$-unavoidable.

This upper bound on the unavoidability of extension of forests is certainly not tight.
In Proposition~\ref{prop:D0n}, we show that the unavoidability of some $k$-extensions of the arcless digraph on $n$ vertices is $2^k n-o(n)$ as $n \to +\infty$. Hence we are left with the following analogue to Problem~\ref{pb:blowup}.

\begin{problem}\label{pb:extension}
    What is the minimum function $f$ such that every $k$-extension of every forest of order $n$ is $(2^{f(k)} \cdot n)$-unavoidable?
\end{problem}
The above-mentioned results show that the function $f$ is at least linear and at most quadratic.
The same question applies to any subfamily, starting with that of empty digraphs.
Generally, determining the precise value or good approximations of the unavoidability of $k$-extensions of oriented trees would be interesting.
In Section~\ref{sec:particular}, we show an almost tight upper bound on the unavoidability of full $k$-twin extensions
and the exact unavoidability of the full $1$-extensions of directed paths.

\section{Preliminaries}

\subsection{Basic properties of tournaments}

We will make use of the following folklore properties of tournaments.

\begin{lemma}[{\cite[Proposition 2.2.2]{bang-jensen2018}}]
\label{prop:degre<k}
Let $k$ be a positive integer.
Every tournament has at most $2k-1$ vertices of out-degree less than $k$.
\end{lemma}

The {\bf out-section} of a vertex $v$ in a tournament $T$, denoted by $S^+_T(v)$  or simply
$S^+(v)$, is the set of vertices that can be reached from $v$
by a directed path.
Similarly, the {\bf in-section} of a vertex $v$ in a tournament $T$, denoted by $S^-_T(v)$ or simply
$S^-(v)$ is the set of vertices from which $v$ can be reached by a directed path.
We set $s^+_T(v)=|S^+_T(v)|$ and $s^-_T(v)=|S^-_T(v)|$.
An {\bf out-generator} (resp. {\bf in-generator}) of $T$ is a vertex
$v$ such that $s^+_T(v) = |V(T)|$ (resp. $s^-_T(v) = |V(T)|$).

\begin{lemma}[{\cite[Proposition~2.7.6]{bang-jensen2018}}]
\label{lem:origin}
Let $T$ be a tournament and $v$ a vertex of $T$.
In $T$, $v$ is the initial (resp. terminal) vertex of a directed path of order
$s^+(v)$ (resp. $s^-(v)$).
In particular, if $v$ is an out-generator (resp. in-generator) of $T$, then $v$
is the initial vertex (resp. terminal vertex) of a directed Hamiltonian path.
\end{lemma}

Let $T$ be a tournament, and $A$ and $B$ be two disjoint sets of vertices of $T$. If every vertex of $A$ dominates every vertex of $B$, then we write $A\Ra B$.

\begin{lemma}\label{lem:reduc-origin}
Let $T$ be a tournament and
let $O$ be the set of vertices that are initial vertices of a directed path of order $k$ in $T$.
Then $O\Ra V(T)\setminus O$.
\end{lemma}

\begin{proof}
    Suppose for contradiction that there exists $u \in V(T) \setminus O$ and $v \in O$ such that $uv \in A(T)$.
    Since $v \in O$, there is a path $P$ of order $k$ in $T$ with initial vertex $v$.
    Hence $s^+(u) \geq |V(P)| \geq k$ and so $u$ is the initial vertex of a directed path of order $k$ by Lemma~\ref{lem:origin}.
\end{proof}

\subsection{Median order}

The notion of median order has been introduced in~\cite{havetJGT35}, see also~\cite[Chapter 2]{bang-jensen2018} and~\cite{draganicCPC30}.
A {\bf median order} of a digraph $D$ is a linear order
$(v_1,v_2, \ldots, v_n)$ of its vertex set such that $|\{(v_i,v_j) \mid i < j\}|$, the number of {\bf forward arcs} ({\it i.e.} directed from left to right), is as large as possible. 
Let us first note some basic properties of median orders of
tournaments.

\begin{lemma}\label{lem:median}
 Let $T$ be a tournament and $(v_1,v_2, \ldots, v_n)$ a
median order of $T$. Then, for any two indices $i,j$ with $1 \leq i < j \leq n$:
\medskip
\begin{enumerate}[label=(M\arabic*)]
\item\label{item:M1} the suborder $(v_i,v_{i+1},\ldots,v_j)$ is a median order of the
  induced subtournament $T\ind{\{v_i,v_{i+1},\ldots,v_j\}}$;
\item\label{item:M2} $v_i$ dominates at least half of the vertices
  $v_{i+1},v_{i+2},\ldots,v_j$, and $v_j$ is dominated by at least half of the vertices $v_i,v_{i+1},\ldots,v_{j-1}$,
  in particular $v_i$ dominates $v_{i+1}$;
  \item\label{item:M3}  if $v_{i+2}\ra v_i$, then $v_{i-1} \Ra \{v_i, v_{i+1}, v_{i+2}\}$ when $i \geq 2$, and $\{v_i, v_{i+1}, v_{i+2}\}\Ra v_{i+3}$ when $i \leq n-3$;
  \item\label{item:M4} if $v_n\ra v_1$, then there exists $\ell$ such that $v_1\ra v_\ell \ra v_n$;
  moreover, if $|\{\ell \mid  v_1\ra v_\ell \ra v_n\}| =1$, then both $(v_2, \ldots, v_n, v_1)$ and $(v_n,v_1, \ldots, v_{n-1})$ are median orders of $T$;
  \item\label{item:M5} if $(v_{i}',\dots,v_{j}')$ is a median order of $T\ind{\{v_i,v_{i+1},\ldots,v_j\}}$, then $(v_1,\ldots, v_{i-1}, v_i',\dots,v_j',v_{j+1},\dots, v_n)$ is also a median order of $T$.
\end{enumerate}
\end{lemma}

\begin{proof}
    Items~\ref{item:M1},~\ref{item:M2}, and~\ref{item:M5} follow easily from the definition.

    \medskip

    \noindent{\it \ref{item:M3}} Assume for a contradiction that $v_{i+2} \to v_i$ but there is an arc from $v_i,v_{i+1},v_{i+2}$ to $v_{i-1}$. 
    By~\ref{item:M2} there is exactly one such arc and its tail $u$ belongs to $\{v_{i+1},v_{i+2}\}$. 
    If $u=v_{i+2}$, then $v_{i+2}$ has two out-neighbours in $\{v_{i-1},v_i,v_{i+1}\}$, contradicting~\ref{item:M2}.
    If $u=v_{i+1}$, then $(v_{i+1},v_{i-1},v_{i+2},v_i)$ has less forward arcs than $(v_{i-1},v_i,v_{i+1},v_{i+2})$, a contradiction to~\ref{item:M1}. 
    By directional duality, we also have $\{v_i, v_{i+1}, v_{i+2}\}\Ra v_{i+3}$.

    \medskip
    
    \noindent{\it \ref{item:M4}}
    Suppose that $v_n \to v_1$.
    By~\ref{item:M2}, $v_1$ dominates at least half of $v_2, \dots, v_n$, and so more than half of $v_2, \dots, v_{n-1}$.
    Symmetrically, $v_n$ is dominated by more than half of $v_2, \dots, v_{n-1}$.
    Hence there exists $\ell \in \{2, \dots, n-1\}$ such that $v_1 \to v_\ell \to v_n$.
    Now suppose that this $\ell$ is unique.
    Hence $v_1$ dominates exactly $\frac{n-1}{2}$ vertices among $v_2, \dots,v_{n}$.
    It follows that the number of forward arcs in $(v_2, \dots, v_n, v_1)$ is 
    $|\{(v_{i},v_{j}) \mid 1 \leq i < j \leq n\}| - d^+(v_1) + d^-(v_1) = |\{(v_{i},v_{j}) \mid 1 \leq i < j \leq n\}|$,
    which is maximum since $(v_1, \dots, v_n)$ is a median order.
    This proves that $(v_1, \dots, v_n,v_1)$ is a median order.
    Symmetrically, $(v_n,v_1, \dots, v_{n-1})$ is a median order too.
\end{proof}

\subsection{A K{\H{o}}v{\'a}ri-S{\'o}s-Tur{\'a}n-like lemma}

The following lemma is inspired by the classical result by K{\H{o}}v{\'a}ri, S{\'o}s, and Tur{\'a}n~\cite{kovariCM3}, and a similar one was recently used in~\cite{draganicCPC30}.

\begin{lemma}\label{lemma:kovaru_sos_turan_adapte}
Let $\alpha, \epsilon$ be two positive reals and $k$ a positive integer.
Let $(A,B,E)$ be a bipartite graph with $|A| \geq \frac{k}{\frac{1}{2}-\epsilon}$ and  $\displaystyle \alpha \leq \frac{\binom{(\frac{1}{2}-\epsilon)|A|}{k}}{\binom{|A|}{k}}$, such that $d(v) \geq (\frac{1}{2}-\epsilon)|B|$ for every $v \in A$. 
Then there exist $A^\star \subseteq A$ and  $B^\star \subseteq B$ such that
\begin{enumerate}
    \item $|A^\star| \geq k$,
    \item $|B^\star| \geq \alpha |B|$, and
    \item $A^\star$ is complete to $B^\star$.
\end{enumerate}
\end{lemma}

\begin{proof}

Suppose for contradiction that the result does not hold.
We will count by two different methods the number $M$ of pairs $(X,b)$ with $X \subseteq A$ of size $k$ and $b \in B$ such that $X \subseteq N(b)$.
On the first hand, for every  $X \subseteq A$ of size $k$, there are  less than $\alpha|B|$ vertices $b$ in $B$ such that 
$X \subseteq N(b)$. We thus have 
\[
M < \binom{|A|}{k}\alpha|B|. 
\]
On the other hand,  by Jensen's inequality, we have
\begin{eqnarray*}
M = \sum_{b \in B}\binom{d(b)}{k} & \geq &  |B| \binom{\frac{1}{|B|}\sum_{b\in B}d(b)}{k} = |B| \binom{\frac{1}{|B|}\sum_{a \in A}d(a)}{k}\\
& \geq & |B| \binom{\frac{|A|}{|B|}(\frac{1}{2} - \epsilon)|B|}{k} \geq \binom{|A|}{k} \alpha |B|.
\end{eqnarray*}
This contradiction proves the lemma.
\end{proof}

Note that $\alpha = 2^{-\frac{k}{\frac{1}{2} - \epsilon}}$ satisfies the hypothesis of Lemma~\ref{lemma:kovaru_sos_turan_adapte}.
Indeed, simple computations show that $\frac{\binom{(\frac{1}{2}-\epsilon)|A|}{k}}{\binom{|A|}{k}}$ is a non-increasing function of $|A|$,
and as $|A| \geq \frac{k}{\frac{1}{2}-\epsilon}$ we have 
\[
  \frac{ \binom{(\frac{1}{2}-\epsilon)|A|}{k}}{  \binom{|A|}{k}} 
  \geq  \frac{ \binom{k}{k}}{\displaystyle \binom{\frac{k}{\frac{1}{2}-\epsilon}}{\scriptstyle k}} 
  \geq 2^{-\frac{k}{\frac{1}{2}-\epsilon}}.
\]

\section{Blow-ups of oriented trees}\label{sec:tree-blow-up}

In this section, we show that for every fixed integer $k$, the class of $k$-blow-ups of oriented trees is linearly avoidable.

\begin{theorem}\label{thm:blowup}
Let $F$ be an oriented tree of order $n$ and $k$ a positive integer.
Every tournament on at least $2^{10+18k} k  n$ vertices contains $F[k]$.
\end{theorem}

\begin{proof}
    We root $F$  on an arbitrary vertex $s$.
    Let $q=2^{1+9 k}$, $p=\lceil 2^{5+3^2 k} k \rceil$ and $\epsilon=\frac{1}{2q}$.
    Observe that we have $\epsilon \leq \frac{1}{6}$, $q \geq 2^{1+ \frac{k}{(\frac{1}{2}-\epsilon)^2}}$, and $p \geq q \frac{k}{(\frac{1}{2}-\epsilon)^2}$.
    Let $N = \frac{1}{2}(4qn + 1)p \leq \frac{1}{2} 2^{10+18k} k \cdot n$ (using $4qn+1 \leq 2 \times 4qn$ and $p \leq 2 \times 2^{5+3^2 k}k$) and let $T$ be a tournament of order $2N$.
    We will show that $T$ has a subdigraph isomorphic to $F[k]$.
    Let $v_{-N}, \dots, v_{N-1}$ be a median order of $T$.
    We partition $v_{-N}, \dots, v_{N-1}$ into $4q n+1$ intervals $V_i=\{v_{pi}, \dots, v_{p(i+1)-1}\}$ for $i\in\{-2q n, \dots, 2q n\}$.
    We will find a copy of $F[k]$ in $D$ iteratively, by adding at each step the out-children (\textit{i.e.} its children which are out-neighbours) or the in-children (\textit{i.e.} its children which are in-neighbours) of a node in $F$.
    Let $F_1, \dots, F_r$ be a sequence of subtrees of $F$ such that
    \begin{enumerate}[label=$(\roman*)$]
        \item $F_1$ is the singleton $\{s\}$,
        \item for every $i \leq \frac{r-1}{2}$, $F_{2i+1}$ is obtained from $F_{2i}$ by adding all the out-children in $F$ of a vertex $x \in V(F_{2i})$ without any out-children in $F_{2i}$,
        \item for every $i \leq \frac{r-2}{2}$, $F_{2i+2}$ is obtained from $F_{2i+1}$ by adding all the in-children in $F$ of a vertex $x \in V(F_{2i+1})$ without any in-children in $F_{2i+1}$, and
        \item $F_r=F$.
    \end{enumerate}
    Observe that such a sequence exists with $r \leq 2n - 1$.
    Moreover, we say that a vertex $x$ in $F_i$ is \textbf{fully processed} if all its children in $F$ are also in $F_i$, \textbf{half processed} if all its out-children in $F$ are in $F_i$. Otherwise we say that $x$ is \textbf{not processed}.
    We prove the following property by induction:
    \begin{center}
        \it
        \begin{minipage}{0.9 \textwidth}
            For every $j \in\{1,\dots, r\}$, there is an injection $\iota:V(F_j) \to \{-2q |V(F_{j})|, \dots, 2q |V(F_{j})|\}$ and for every $x \in V(F_j)$
            there exists a set $U_x \subseteq V_{\iota(x)}$ such that
            \begin{enumerate}[label=$(\roman*)$]
                \item $|U_x| \geq \frac{k}{\left(\frac{1}{2}-\epsilon\right)^2}$ if $x$ is not processed,
                \item $|U_x| \geq \frac{k}{\frac{1}{2}-\epsilon}$ if $x$ is half processed,
                \item $|U_x| \geq k$ if $x$ is fully processed,
                \item if $xy$ is an arc in $F_j$, then every vertex in $U_x$ dominates $U_y$, and
                \item for every $i \in \{-2q |V(F_{j})|, \dots, 2q |V(F_{j})|\}$ a fraction of at most $\frac{1}{q}$ of the elements of the
            interval $\{i, \dots, 2q |V(F_{j})|\}$ (respectively $\{-2q |V(F_{j})|, \dots, i\}$) is in $\{\iota(x) \mid x \in V(F_j)\}$.
            More precisely, 
            \begin{align*}
                |\{i,\dots, 2q |V(F_{j})|\} \cap \{\iota(x) \mid x \in V(F_j)\}| &\leq \frac{1}{q} \left(2q |V(F_{j})|-i+1\right),\\
                \text{and~~~~} |\{-2q |V(F_{j})|,\dots,i\} \cap \{\iota(x) \mid x \in V(F_j)\}| &\leq \frac{1}{q} \left(i+2q |V(F_{j})|+1\right).
            \end{align*}
            \end{enumerate}
        \end{minipage}
    \end{center}

    Observe that, for $j=1$, setting $\iota(s)=0$ and taking for $U_s$ any set in $V_0$ of size $p \geq \frac{k}{\left(\frac{1}{2}-\epsilon\right)^2}$ satisfies the invariant. Note also that, for $j=r$, the invariant implies the theorem.
    Now suppose that the property holds for some $j<r$.
    Let $\iota$ and $(U_x)_{x \in F_j}$ witnessing that the property holds at step $j$.
    Let $x$ be the vertex in $F_j$ such that $F_{j+1}$ is obtained from $F_j$ by adding the out-children of $x$ in $F$, or the in-children of $x$ in $F$.
    Our goal is to find a new value $U'_x$ for $U_{x}$, and values for $\iota(y),U_y$ for every out-children (respectively in-children) $y$ of $x$. 

    \begin{description}
        \item[{\bf Case 1:}] {\bf $F_{j+1}$ is obtained from $F_j$ by adding the $\ell$ out-children of $x$ in $F$.}
        
    Let $A = U_{x}$  and $B = \bigcup_{\iota(x) < i \leq 2q |V(F_{j})|+q\ell} V_i$.
    Note that $|B| \geq p(2q |V(F_{j})| + \ell q - \iota(x)) \geq \ell q p$. 
    We first justify that every vertex in $V_{\iota(x)}$ dominates at least $(\frac{1}{2} - \epsilon)|B|$ vertices in $B$. To this purpose, let $v_i$ be any vertex in $V_{\iota(x)}$, so $i\in \{p\iota(x),\dots,p(\iota(x)+1)-1\}$. Let $X$ be the set of vertices $\{v_{i+1},\dots,v_{p(\iota(x)+1)-1}\}$, which is disjoint from $B$ ($X$ is empty if $i=p(\iota(x)+1)-1$). Observe that
    \[ 
    N^+(v_i) \cap B = \Big(N^+(v_i) \cap (X \cup B)\Big) \setminus \Big(N^+(v_i) \cap X \Big).
    \]
    Since $X\cup B = \{v_{i+1},\dots,v_{p(\iota(x)+1)-1+|B|}\}$, and because $v_{-N},\dots,v_{N-1}$ is a median order, by \ref{item:M2} we have $|N^+(v_i) \cap (X \cup B)| \geq \frac{1}{2}|X \cup B|$.
    Therefore, we obtain
    \[ 
    |N^+(v_i) \cap B| \geq \frac{1}{2}|X \cup B| - |X| = \frac{1}{2}|B| - \frac{1}{2} |X| \geq \frac{1}{2}(|B|-p) \geq \left(\frac{1}{2}-\frac{1}{2q}\right)|B| = \left(\frac{1}{2}-\epsilon \right)|B|.
    \]
    This shows that every vertex in $V_{\iota(x)}$ dominates at least $(\frac{1}{2} - \epsilon)|B|$ vertices in $B$, and so does every vertex in $A$. By Lemma~\ref{lemma:kovaru_sos_turan_adapte} applied for $\alpha = 2^{-\frac{k}{(\frac{1}{2} - \epsilon)^2}}$, 
    there exist $A^\star \subseteq U_{x}$ and $B^\star \subseteq B$ such that every vertex in $A^\star$ dominates $B^\star$, 
    $|B^\star| \geq 2^{-\frac{k}{(\frac{1}{2}-\epsilon)^2}}|B|$, and $|A^\star| \geq \frac{k}{\frac{1}{2}-\epsilon}$.
    We set $U'_x = A^\star$.
    
    Let $I$ be the set of indices $i\in \{\iota(x)+1, \dots, 2q |V(F_{j})| + q \ell\}$ such that $|B^{\star}\cap V_i| \geq \frac{k}{(\frac{1}{2} - \epsilon)^2}$ and $i \notin \{\iota(z) \mid z \in V(F_j)\}$. In what follows we show that $I$ is large enough, which will enable us to find $U_y$ for every out-child $y$ of $x$ in $F$.

    \begin{claim}\label{claim:find_emplacement_for_children_of_x}
        $|I| \geq \ell$
    \end{claim}
    \begin{proofclaim}
        Suppose that this does not hold, so $|I| \leq \ell-1$. We define $J, I^\star,$ and $I^\iota$ as follows:
        \begin{align*}
            J &= \{\iota(x)+1, \dots, 2q |V(F_{j})| + ql\},\\
            I^\star &= \{ i\in J \mid |V_i \cap B^\star | \geq \tfrac{k}{(\frac{1}{2}-\epsilon)^2}\},\text{ and}\\
            I^\iota &= \{ \iota(z) \mid z\in V(j) \} \cap J.
        \end{align*}
        Observe that $I = I^\star \setminus I^\iota$ by definition. 
        By the induction hypothesis, $|I^\iota| \leq \frac{1}{q}(|J|-q\ell)$. Since $|I| \leq \ell-1$ by assumption, we obtain $|I^\star| \leq \ell-1 + \frac{1}{q}(|J|-q\ell) = \frac{|J|}{q} - 1$. This leads to the following inequalities, in which we also use $|B| = p|J|$
        \[
        \begin{split}
            |B^\star| &= \sum_{i\in J\cap I^\star}|V_i \cap B^\star| + \sum_{i\in J\setminus I^\star}|V_i \cap B^\star|\\
            &\leq \left(\frac{|J|}{q}-1\right)p + \left( \frac{q-1}{q}|J|+1\right)\frac{k}{(\frac{1}{2}-\epsilon)^2} \\
            &= \left(\frac{|B|}{q}\right) + \left(\frac{q-1}{q}\frac{|B|}{p}\frac{k}{(\frac{1}{2}-\epsilon)^2}\right) + \left(\frac{k}{(\frac{1}{2}-\epsilon)^2} - p\right)\\
            &< \left( \frac{1}{q} + \frac{k}{p(\frac{1}{2}-\epsilon)^2} \right)|B|\\
            &\leq \frac{2}{q} |B| \\
            &\leq 2^{-\frac{k}{(\frac{1}{2}-\epsilon)^2}}|B|.
        \end{split}
        \]
        This is a contradiction since $|B^\star| \geq 2^{-\frac{k}{(\frac{1}{2}-\epsilon)^2}}|B|$.
    \end{proofclaim}

    We can now extend $\iota$ to the out-children of $x$ by taking values in $I$, and then we can arbitrarily take $U_z \subseteq B^\star \cap V_{\iota(z)}$
    for every out-child $z$ of $x$.
    Recall that $|V(F_{j+1})| = |V(F_{j})| + \ell$. The invariant still holds because we add $2\ell q p$ vertices at the end and at the beginning of the median order:

    \begin{enumerate}[label=$(\roman*)$]
        \item For every $i \in \{-2q |V(F_{j})|-2q\ell, \dots, 2q |V(F_{j})|\}$, $|\{i,\dots, 2q |V(F_{j})|+2q \ell\} \cap \{\iota(z) \mid z \in V(F_{j+1})\}|
                         \leq \frac{1}{q} (2q |V(F_{j})|-i+1) + \ell \leq \frac{1}{q}(2q |V(F_{j+1})|-i+1)$.
        \item For every $i \in \{2q |V(F_{j})|+1, \dots, 2q |V(F_{j})|+q \ell\}$, $|\{i,\dots, 2q |V(F_{j})|+2q \ell\} \cap \{\iota(z) \mid z \in V(F_{j+1})\}| 
                         \leq \ell \leq \frac{1}{q}(q \ell) \leq \frac{1}{q}(2 q |V(F_{j+1})|-i+1)$.
        \item For every $i \in \{2q |V(F_{j})|+q\ell+1, \dots, 2q |V(F_{j})|+2q \ell\}$, 
                        $|\{i,\dots, 2q |V(F_{j})|+2q \ell\} \cap \{\iota(z) \mid z \in V(F_{j+1})\}| = 0 \leq \frac{1}{q}(2 q |V(F_{j+1})|-i+1)$.
        \item For every $i \in \{-2q |V(F_{j})|, \dots, 2q |V(F_{j})|+2q\ell\}$, $|\{-2q |V(F_{j})|-2q \ell, \dots, i\} \cap \{\iota(z) \mid z \in V(F_{j+1})\}|
                         \leq \frac{1}{q} (i+2q |V(F_{j})|+1) + \ell \leq \frac{1}{q}(i+2q |V(F_{j+1})|+1)$.
        \item For every $i \in \{-2q \ell - 2q |V(F_{j})|, \dots, -2q |V(F_{j})|-1\}$, $|\{-2q |V(F_{j})|-2q \ell, \dots, i\} \cap \{\iota(z) \mid z \in V(F_{j+1})\}| = 0
                        \leq \frac{1}{q}(i+2q |V(F_{j+1})|+1)$.
    \end{enumerate}

    \item[{\bf Case 2:}] {\bf$F_{j+1}$ is obtained from $F_j$ by adding the $\ell$ in-children of $x$ in $F$.}
    
    Let $A = U_{\iota(x)}$ and $B = \bigcup_{-2q |V(F_{j})|-q\ell \leq i < \iota(x)} V_i$.
    Since $x$ is half processed, we have $|A| \geq \frac{k}{\frac{1}{2}-\epsilon}$.
    Moreover $|B| \geq q \ell p$ and every vertex in $A$ is dominated by at least 
    $\frac{1}{2}|B| - |A| \geq (\frac{1}{2}-\epsilon)|B|$ vertices in $B$.
    Using Lemma~\ref{lemma:kovaru_sos_turan_adapte},
    we deduce that there are sets $A^\star\subseteq A, B^\star \subseteq B$ with $B^\star \Rightarrow A^\star$, $|A^\star| \geq k$ and $|B^\star| \geq 2^{-\frac{k}{\frac{1}{2}-\epsilon}}|B| \geq 2^{-\frac{k}{(\frac{1}{2}-\epsilon)^2}}|B|$.
    We set $U'_x = A^\star$.
    We now find $U_y$ for every $y$ in-children of $x$ in $F$.
    With a proof almost identical to the one of Claim~\ref{claim:find_emplacement_for_children_of_x} we deduce the following claim.
    \begin{claim}
        There is a set $I \subseteq \{-2q |V(F_{j})| - q \ell, \dots, \iota(x)-1\}$ of size at least $\ell$ such that
        $\iota(z) \not\in I$ for every $z \in V(F_j)$ and $|B^\star \cap V_i| \geq \frac{k}{(\frac{1}{2}-\epsilon)^2}$ for every $i \in I$.
    \end{claim}
    Then we build $U_y, \iota(y)$ for every in-children $y$ of $x$ as in the first case.
    This concludes the proof of the theorem.\qedhere
    \end{description}
\end{proof}

\section{Extensions of oriented trees}\label{sec:extension-tree}

In this section, we show that for every fixed integer $k$, the class of $k$-extensions of oriented trees is linearly unavoidable.
Our method is based on the fact that the number of copies of $TT_{2k+1}$ in any $N$-vertex tournaments is in $\Omega_k(N^{2k+1})$.
We start with the case $k=1$, which is a good warm-up for the general case.

\begin{lemma}[Folklore]\label{lemma:count_TT3}
Every tournament of order $N$ contains at least $\frac{1}{8} N(N-1)(N-3)$ copies of $TT_3$.
\end{lemma}

\begin{proof}
Let $T$ be a tournament of order $N$.
    Let $\vec{t}$ be the number of directed triangles in $T$ and let $tt_3$ be the number of $TT_3$ in $T$.
    The number of copies of $\vec{P}_3$ (the directed path of order $3$) in $T$ is $\sum_{v \in V(T)} d^-(v)d^+(v) \leq N \left(\frac{N-1}{2}\right)^2 = \frac{1}{4} N(N-1)^2$.
    But the number of $\vec{P}_3$ is also $tt_3 + 3\vec{t} = tt_3 + 3(\binom{N}{3}-tt_3) = 3\binom{N}{3} - 2tt_3$.
    We deduce that $3\binom{N}{3} - 2tt_3 \leq \frac{1}{4} N(N-1)^2$ and so $tt_3 \geq \frac{1}{2}\left(\frac{N(N-1)(N-2)}{2}-\frac{N(N-1)^2}{4}\right)= \frac{1}{8} N(N-1)(N-3)$.
\end{proof}

We will also need the following well-known lemma.
Recall that a \emph{hypergraph} $\mathcal{H}$ is a pair $(V,E)$ where $V$ is a finite set and $E$ is a subset of $2^{V}$.
We call the elements of $V$ the \emph{vertices} of $\mathcal{H}$, and the elements of $E$ the \emph{hyperedges} of $\mathcal{H}$.
An hypergraph $\mathcal{H}' = (V',E')$ is a \emph{subhypergraph} of $\mathcal{H}$, and we write $\mathcal{H}' \subseteq \mathcal{H}$, if $V' \subseteq V$ and $E' \subseteq E$.

\begin{lemma}[Folklore]\label{lemma:from_avg_degree_to_min_degree}
    Let $\mathcal{H}$ be a hypergraph on $n$ vertices and with $m$ hyperedges.
    There is a hypergraph $\mathcal{H}' \subseteq \mathcal{H}$ such that every vertex of $\mathcal{H}'$ belongs to at least $\frac{m}{n}$ hyperedges of $\mathcal{H}'$.
\end{lemma}

\begin{proof}
    We proceed by induction on $n$. If every vertex in $\mathcal{H}$ is in at least $\frac{m}{n}$ hyperedges we are done.
    Otherwise there is a vertex $u$ in less than $\frac{m}{n}$ hyperedges, then $\mathcal{H} - u$ has $m' \geq m-\frac{m}{n} = \frac{m}{n} (n-1)$ hyperedges and $n'=n-1$ vertices.
    Hence, by the induction hypothesis, $\mathcal{H}-u$ has a subhypergraph with minimum degree at least $\frac{m'}{n'} \geq \frac{m}{n}$ as wanted.
\end{proof}

Given an acyclic digraph $D$, we denote by $\unvd^* (D)$ the maximum of $\unvd(H) + |V(D) \setminus V(H)|$ over all the subdigraphs $H$ of $D$.

\begin{theorem}\label{thm:1-extension_of_trees}
    Let $D$ be a $1$-extension of an oriented tree $F$ on $n$ vertices, then 
    $\unvd(D) \leq 8\unvd^*(F) + 3$.
\end{theorem}

By applying Theorem~\ref{thm:kuhn_mycroft_osthus} to bound $\unvd^*(F)$, we obtain the following.

\begin{corollary}
    There exists a constant $n_0$ such that the following holds.
    Let $D$ be a $1$-extension of an oriented tree $F$ on at least $n_0$ vertices, then $\unvd(D) \leq 16|V(F)|-13$.
\end{corollary}

An {\bf embedding} $\phi\colon H \to D$ of a digraph $H$ in a digraph $D$ is an injective function $\phi\colon V(H) \to V(D)$ such that
$\phi(u)\phi(v) \in A(D)$ for every $uv \in A(H)$.
Observe that $G$ contains a copy of $H$ if and only if there is an embedding of $H$ in $G$.

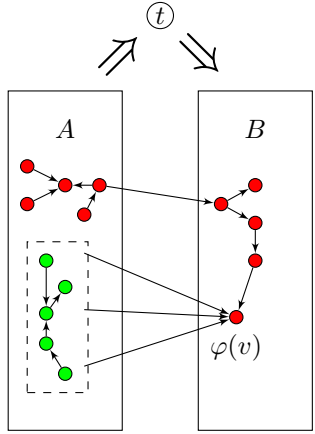
\begin{figure}[hbt!]
    \centering
    \begin{tikzpicture}
        \tikzset{vertex/.style = {circle,draw,minimum size=4pt, inner sep=1pt}}
        \tikzset{vertexblack/.style = {circle,draw,fill=red,minimum size=5pt, inner sep=0pt}}
        \tikzset{edge/.style = {->,> = latex'}}

        \node[vertex] (v) at (0,0) {$t$};

        \draw (-2,-5.5) rectangle (-0.5,-1);
        \draw (0.5,-5.5) rectangle (2,-1);

        \node at (-1.25, -1.5) {$A$};
        \node at (1.25, -1.5)  {$B$};

        \node[vertexblack,label=below:$\phi(v)$] (x) at (1,-4) {};
        \draw[edge] (-1,-3.15) to (x);
        \draw[edge] (-1,-3.9) to (x);
        \draw[edge] (-1,-4.65) to (x);

        \node[rotate=45] at (-0.5,-0.55) {\huge $\Rightarrow$};
        \node[rotate=-45] at (0.5,-0.55) {\huge $\Rightarrow$};

        \node[vertexblack] (x1) at (-1.75, -2) {};
        \node[vertexblack] (x2) at (-1.25, -2.25) {};
        \node[vertexblack] (x3) at (-1.75, -2.5) {};
        \node[vertexblack] (x4) at (-0.8, -2.25) {};
        \node[vertexblack] (x5) at (-1, -2.64) {};
        \draw[edge] (x1) to (x2);
        \draw[edge] (x3) to (x2);
        \draw[edge] (x5) to (x4);
        \draw[edge] (x4) to (x2);

        \node[vertexblack] (x6) at (0.8, -2.5) {};
        \node[vertexblack] (x7) at (1.25, -2.25) {};
        \node[vertexblack] (x8) at (1.25, -2.75) {};
        \node[vertexblack] (x9) at (1.25, -3.25) {};
        \draw[edge] (x6) to (x7);
        \draw[edge] (x6) to (x8);
        \draw[edge] (x8) to (x9);
        \draw[edge] (x9) to (x);

        \draw[edge] (x4) to (x6);

        \draw[dashed] (-1.75,-3) rectangle (-0.95,-5);

        \node[vertexblack, fill=green] (z1) at (-1.5, -3.25) {};
        \node[vertexblack, fill=green] (z2) at (-1.25, -3.6) {};
        \node[vertexblack, fill=green] (z3) at (-1.5, -3.95) {};
        \node[vertexblack, fill=green] (z4) at (-1.5, -4.35) {};
        \node[vertexblack, fill=green] (z5) at (-1.25, -4.75) {};

        \draw[edge] (z3) to (z2);
        \draw[edge] (z1) to (z3);
        \draw[edge] (z4) to (z3);
        \draw[edge] (z5) to (z4);
    \end{tikzpicture}
    \caption{Illustration of the proof of Theorem~\ref{thm:1-extension_of_trees}. The red vertices are the image of $V(F-Z)$ under $\phi$. The green vertices represent an embedding of $F\ind{Z}$ in the set $A \cap N^-_T(\phi(v)) \setminus \phi(V(F-Z))$, which is dashed in the figure.}
    \label{fig:1extension_of_a_tree}
\end{figure}

\begin{proof}[Proof of Theorem~\ref{thm:1-extension_of_trees}]
    Let $D$ be a $1$-extension of $F$ with added vertex $s$.
    We assume that $D$ is a full $1$-extension of $F$, for otherwise we just add missing arcs.
    Let $T$ be a tournament on $N=8\unvd^*(F)+3$ vertices.
    By Lemma~\ref{lemma:count_TT3}, the number of $TT_3$ in $T$ is at least $\frac{1}{8}N(N-1)(N-3)$.
    We deduce that there is a vertex $t$ which is the middle vertex (\textit{i.e.} the vertex with in- and out-degree $1$) of at least $\frac{1}{8}(N-1)(N-3)$ copies of $TT_3$ in $T$.
    In other words there are at least $\frac{1}{8}(N-1)(N-3)$ arcs from $N^-(t)$ to $N^+(t)$.
    Consider the bipartite graph $H_0=(N^-_T(t),N^+_T(t),E)$ where $E = \{uw \in A(T) \mid u \in N^-_T(t), w \in N^+_T(t)\}$.
    By Lemma~\ref{lemma:from_avg_degree_to_min_degree}, $H_0$ has a nonempty subgraph $H$ with minimum degree at least $\frac{(N-1)(N-3)}{8(N-1)} = \frac{1}{8}(N-3) \geq \unvd^*(F)$.
    Let $A = V(H) \cap N^-_T(v), B = V(H) \cap N^+_T(t)$. Since $H$ is bipartite, note that in particular $|A|,|B| \geq \unvd^*(F)$.

    We will now inductively build and embedding $\phi\colon F \to T\ind{V(H)}$ such that
    \begin{itemize}
        \item for every $y \in N^-_D(s)$, $\phi(y) \in A$, and
        \item for every $y \in N^+_D(s)$, $\phi(y) \in B$.
    \end{itemize}
    Extending such an embedding of $F$ by setting $\phi(s) = t$ yields the desired embedding of $D$ in $T$.
    
    If $F$ is empty, then the result is clear.
    Suppose now that $F$ is not empty and consider the tree $F'$ obtained from $F$ by contracting every connected component $C$ of $D\ind{N^-_D(s)}$ and $D\ind{N^+_D(s)}$ into a single vertex $u_{V(C)}$.
    Consider a leaf $u_Z$ of $F'$, and without loss of generality, assume that $Z \subseteq N^-_D(s)$.
    If $Z=V(F)$, then $V(F) \subseteq N^-_D(s)$, and we can embed $F$ in $T\ind{A}$ since $|A| \geq \unvd^*(F) \geq \unvd(F)$.
    Otherwise, assume by induction that we have an embedding $\phi\colon F-Z \to T\ind{V(H)}$ with $\phi(N^-_D(s) \setminus Z) \subseteq A$ and $\phi(N^+_D(s)) \subseteq B$.
    Let $uv$ be the arc joining $Z$ to $F-Z$ in $F'$. 
    Observe that we have $u \in Z$ and $v \in N^+_D(s)$, for otherwise $uvs$ is a directed triangle in $D$, a contradiction to $D$ being acyclic.
    
    As $H$ has minimum degree at least $\unvd^*(F)$, and because $\phi(v) \in B$, the size of $A \cap N^-_T(\phi(v))$ is at least $\unvd^*(F)$, hence 
    there are at least $\unvd^*(F)-|V(F) \setminus Z| \geq \unvd(F\ind{Z})$ vertices in $A \cap N^-_T(\phi(v)) \setminus \phi(V(F-Z))$.
    Thus we can extend $\phi$ to the vertices in $F\ind{Z}$ in $A \cap N^-_T(\phi(v)) \setminus \phi(V(F-Z))$, and this gives the desired embedding of $D$.
    See Figure~\ref{fig:1extension_of_a_tree} for an illustration.
\end{proof}

We now move to the general case. We first prove an analogue of Lemma~\ref{lemma:count_TT3} for larger values of $k$.

\begin{lemma}\label{lemma:number_of_TT}
    Let $k$ be a positive integer and let $T$ be a tournament on $N \geq 3^{\binom{k+1}{2}}$ vertices. Then the number of distinct copies of $TT_k$ in $T$ is at least
    $3^{-\binom{k+1}{2}}N^k$.
\end{lemma}

\begin{proof}
    We proceed by induction on $k$.
    For $k=1$ the result is clear.
    Now suppose $k>1$ and that the result holds for smaller values of $k$.
    By Proposition~\ref{prop:degre<k}, in $T$ there are at least $\frac{N}{3}$ vertices of out-degree at least $\frac{N}{3}$.
    For every such vertex $u$, there are by induction at least $3^{-\binom{k}{2}}\left(\frac{N}{3}\right)^{k-1}$ copies of $TT_{k-1}$ in the subtournament induced by $N^+(u)$.
    Together with $u$, this gives as many copies of $TT_k$ in $T$.
    We conclude that the number of copies of $TT_k$ in $T$ is at least $\frac{N}{3} \times 3^{-\binom{k}{2}} \left(\frac{N}{3}\right)^{k-1} = 3^{-\binom{k+1}{2}}N^k$.
\end{proof}

This bound could be improved to $(2^{-\binom{k}{2}}+o(1)) N^k$, see~\cite{coreglianoJGT85}. 

\begin{theorem}\label{thm:k-extension-of-trees}
    Let $k,n$ be integers with $k \geq 1, n \geq 2$ and let $F$ be an oriented tree of order $n$.
    Let $D$ be a $k$-extension of $F$,
    then $\unvd(D) \leq 3^{\binom{2k+2}{2}} \cdot \unvd^*(F)$.
\end{theorem}

Again, applying Theorem~\ref{thm:kuhn_mycroft_osthus} to bound $\unvd^*(F)$, we deduce the following corollary.

\begin{corollary}\label{cor:k-extension-of-trees}
    There is a constant $n_0$ such that the following holds.
    Let $D$ be a $k$-extension of an oriented tree $F$ on at least $n_0$ vertices, then
    $\unvd(D) \leq 2 \cdot 3^{\binom{2k+2}{2}} \cdot |V(F)|$ .
\end{corollary}

\begin{proof}[Proof of Theorem~\ref{thm:k-extension-of-trees}]
    Let $T$ be a tournament on $N=3^{\binom{2k+2}{2}}\unvd^*(F)$ vertices.
    
    Let $H$ be a copy of $TT_{2k+1}$ with acyclic ordering $(u'_0,u_1,u'_1, \dots, u_k,u'_k)$.
    By Lemma~\ref{lemma:number_of_TT}, $T$ contains at least $3^{2k\binom{2k+2}{2}}\unvd^*(F)^{2k+1}$ copies of $H$.
    By the Pigeon-hole Principle, there are $k$ vertices $v_1,\dots, v_k$ in $T$ such that
    there are at least $\frac{1}{N^k} 3^{2k\binom{2k+2}{2}}\unvd^*(F)^{2k+1} = 3^{k\binom{2k+2}{2}}\unvd^*(F)^{k+1}$ embeddings of $H$ in $T$ 
    with $u_i$ mapped to $v_i$ for every $i\in\{1,\dots, k\}$.
    From now on, let us say that an embedding $\phi \colon H \to T$ is {\bf valid} if $\phi(u_i) = v_i$ for every $i\in\{1,\dots, k\}$.
    
    We consider the hypergraph $\mathcal{H}$ with vertex set $V(T) \setminus\{v_1, \dots, v_k\}$ which contains as a hyperedge $\{\phi(u_i') \mid 0\leq i \leq k\}$ for every valid embedding $\phi \colon H \to T$.
    By Lemma~\ref{lemma:from_avg_degree_to_min_degree} applied to $\mathcal{H}$, there is a set of vertices $X \subseteq V(\mathcal{H})$
    such that every vertex in $X$ has degree at least $\frac{1}{N}3^{k\binom{2k+2}{2}}\unvd^*(F)^{k+1}$ in $\mathcal{H}$. 
    Let $(X_0,\dots,X_k)$ be the partition of $X$ where $X_i$ contains all the vertices $x\in X$ 
    such that $x=\phi(u_i')$ for some valid embedding $\phi \colon H\to T$. 
    
    For all $i,j \in \{0, \dots, k\}$ with $i<j$ and for every vertex $x\in X_i$, we claim that $x$ has at least $\frac{1}{N^k} 3^{k\binom{2k+2}{2}}\unvd^*(F)^{k+1} = \unvd^*(F)$ out-neighbours in $X_j$. Assume that this is not the case, then the degree of $x$ in $\mathcal{H}\ind{X}$, which is exactly the number of valid embedding $\phi \colon H\to T$ with the extra conditions $\phi(u_i') = x$ and $\phi(u_\ell') \in X$ for every $\ell\in \{0,\dots,k\}$, is at most 
    \[|N^+(x) \cap X_j| \cdot \prod_{\ell \in \{0, \dots, k\} \setminus \{i,j\}} |X_\ell|  < N^{k-1} \frac{1}{N^k} 3^{k\binom{2k+2}{2}}\unvd^*(F)^{k+1}, \]
    a contradiction to the choice of $X$. Similarly, if $j<i$ then $x$ has at least $\unvd^*(F)$ in-neighbours in $X_j$.
    Since $k\geq 1$, this implies in particular that $|X_i| \geq \unvd^*(F)$ for every $i \in \{0,\dots,k\}$.

    \medskip

    Let $D$ be a full $k$-extension of $F$. As $D$ is acyclic, let $\sigma$ be a topological ordering of $V(D)$. Let $W \subseteq V(D)$ be such that $F = D-W$ and $w_1,\dots,w_k$ the ordering of $W$ with respect to $\sigma$.
    We let $Y_0,\dots,Y_k$ be the partition of $V(F)$ where $Y_i$ contains every vertex between $w_{i}$ and $w_{i+1}$ in $\sigma$,
    for every $i \in \{1,\dots,k-1\}$,
    $Y_0$ contains every vertex before $w_1$ in $\sigma$,
    and $Y_k$ contains every vertex after $w_k$ in $\sigma$.
    By construction, every arc between $Y_i$ and $Y_j$ in $F$ is from $Y_i$ to $Y_j$ if $i<j$.
    We will now construct by induction on $|V(F)|$ an embedding $\phi\colon F \to T\ind{X}$ such that $\phi(Y_i) \subseteq X_i$ for every $i\in\{0,\dots,k\}$.
    By setting $\phi(w_i)=v_i$ for every $i \in [k]$, we will get the desired embedding of $D$ in $T$.

    If $F$ is empty, then the result is clear. Now assume $|V(F)|>0$.
    Consider the tree $F'$ obtained from $F$ by contracting every connected component $C$ of $D\ind{Y_i}$ into a single vertex $x_{V(C)}$, for every $i\in \{0,\dots,k\}$. 
    Let $x_Z$ be a leaf of $F'$ and $i \in \{0,\dots,k\}$ be the index such that $Z \subseteq Y_i$.
    If $Z=V(F)$, then, as $|X_i| \geq \unvd^*(F) \geq \unvd(F)$, there is an embedding $\phi$ of $F$ into $T\ind{X_i}$ and we are done. 
    Otherwise let $p$ be the unique neighbour in $F$ of a vertex in $Z$ which does not belong to $Z$. 
    Assume that $p$ is the tail of the (unique) arc between $F- Z$ and $Z$,
    the other case being symmetric. 
    Let $j \in \{0, \dots,k\}$ be the index such that $p \in Y_j$.
    Recall that, by construction, we have $j<i$.
    By induction hypothesis, there is an embedding $\phi\colon F-Z \to T\ind{X_0\cup \dots \cup X_k}$ such that $\phi(Y_\ell\setminus Z) \subseteq X_\ell$ for every $\ell\in \{0,\dots,k\}$.
    We know that $\phi(p)$ has at least $\unvd^*(F)$ in-neighbours in $X_i$.
    Moreover $\unvd^*(F) - |V(F)\setminus Z| \geq \unvd(F\ind{Z})$ and so $|X_i \cap N^-(\phi(p)) \setminus \phi(V(F-Z))| \geq \unvd(F\ind{Z})$ (by definition of $\unvd^*(F)$).
    It follows that $F\ind{Z}$ can be embedded in $X_i \cap N^-(\phi(p)) \setminus \phi(V(F-Z))$. 
    This gives the desired embedding of $F$ and concludes the proof of the theorem.    
\end{proof}

\section{Particular extensions of some particular forests}\label{sec:particular}

\subsection{Twin extensions of arcless digraphs}

For every positive integers $k,n_1, n_2$, let $D_k(n_1,n_2)$ be full $k$-twin-extension of the arcless digraph $D$ of order $n_1+n_2$ in which the $k$ added vertices have the same in-neighbourhood $V_1$ of size $n_1$ and the same out-neighbourhood $V_2$ of size $n_2$ (with $V_1,V_2$ a partition of $V(D)$).
Note that $D_1(n_1, n_2)$ is simply the oriented star made of a vertex dominated by $n_1$ vertices and dominating $n_2$ vertices.

\begin{proposition}
\label{prop:unvd_star}
For every positive integers $n,n_1,n_2$ with $n=n_1+n_2+1$, each of the following holds
\begin{enumerate}
    \item $\unvd(D_1(0,n-1)) = \unvd(D_1(n-1,0)) = 2n-2$, and
    \item $\unvd(D_1(n_1,n_2)) =  2n -3$.
\end{enumerate}
\end{proposition}

\begin{proof}
    The fact that $\unvd(D_1(0,n-1)) = \unvd(D_1(n-1,0))  \leq 2n-2$ and $\unvd(D_1(n_1,n_2)) \leq  2n -3$ follows from Lemma~\ref{prop:degre<k}.
    For the lower bounds, consider first any $(n-2)$-regular tournament $T$.
    Then $T$ has $2n-3$ vertices and does not contain $D_1(0,n-1)$. Thus $\unvd(D_1(0,n-1)) = \unvd(D_1(n-1,0)) \geq 2n-2$.
    Similarly, if $T_1$ is an $(n_1-1)$-regular tournament and $T_2$ an $(n_2-1)$-regular tournament,
    then $T_1 \Rightarrow T_2$ is a tournament of order $2n-4$ which does not contain $D_1(n_1,n_2)$, implying that $\unvd(D_1(n_1,n_2)) \geq 2n -3$.
\end{proof}

To prove a lower bound on $\unvd(D_k(0,n))$, we will use a probabilistic argument relying on Chernoff's bound.

\begin{proposition}[Chernoff's Bound]\label{chernoff}
    If $X$ is a random variable following a binomial law with parameters $p \in [0,1]$ and $n \geq 0$, then
    for every $\epsilon \in [0,1]$
    \[
    \Pr[X \geq (1+\epsilon)pn] \leq \exp \left(-\frac{\epsilon^2}{3} pn\right).
    \]
\end{proposition}

\begin{proposition}\label{prop:D0n}
    For every fixed $k \geq 1$, $\unvd(D_k(0,n)) = 2^k n - o(n)$ as $n \to +\infty$.
\end{proposition}

\begin{proof}
    The added vertices of $D_k(0,n)$ are sources, so by Proposition~\ref{prop:double}, $\unvd(D_k(0,n)) \leq 2^k n$.
    
    To prove the lower bound, let $ 0 < \epsilon < 1$ and let $T$ be a tournament of order $\left\lfloor \frac{2^k}{1+\epsilon}n+k \right\rfloor$ taken uniformly at random.
    For every $X \subseteq V(T)$ of size $k$, let $d^+_X$ be the number of common out-neighbours of the vertices in $X$, that is $d^+_X = |\bigcap_{x \in X} N^+(x)|$.
    Then $d^+_X$ follows a binomial law with parameters $2^{-k}$ and $\left\lfloor \frac{2^k}{1+\epsilon}n \right\rfloor \leq \frac{2^k}{1+\epsilon} n$.
    By Chernoff's bound (Proposition~\ref{chernoff}), 
    \[
    \begin{split}
      \Pr[d^+_X\geq n] 
      &\leq \Pr\big[d^+_X\geq (1+\epsilon) \E[d^+_X]\big] \\
      &\leq \exp\left( -\frac{\epsilon^2}{3} \E[d^+_X] \right) \\
      &\leq \exp\left( -\frac{\epsilon^2}{3} 2^{-k} \left\lfloor \frac{2^k}{1+\epsilon}n \right\rfloor\right) \\
      &\leq    \exp\left( -\frac{\epsilon^2}{3(1+\epsilon)}n +1\right). \\
    \end{split}
    \]
    It follows by the Union Bound that $\displaystyle \Pr[T\text{ contains }D_k(0,n)] \leq \binom{\lfloor \frac{2^k}{1+\epsilon}n \rfloor}{k}\e^{-\frac{\epsilon^2}{3(1+\epsilon)} n + 1}<1$ for $n$ large enough.
    This proves that $\unvd(D_k(0,n)) > \left\lfloor \frac{2^k}{1+\epsilon}n+k \right\rfloor$ for $n$ large enough.
\end{proof}

We will now prove a similar upper bound for $D_k(n_1,n_2)$.
To do so, we will use Ramsey's theorem, that we first recall here.

\begin{theorem}[Ramsey's theorem~\cite{ramsey1987}]
    For every positive integers $k$, $a$ and $b$, there exists $R_k(a,b)$ such that for every $2$-colouring of the $k$-subsets of $[R_k(a,b)]$,
    either there is a set $X \subseteq [R_k(a,b)]$ of size $a$ such that every $k$-subset of $X$ is coloured $1$, or there is a set
    $Y \subseteq [R_k(a,b)]$ of size $b$ such that every $k$-subset of $Y$ is coloured $2$.
\end{theorem}

\begin{theorem}\label{theorem:k-extended-stable-set}
    Let $\epsilon>0$ and let $k$ be a positive integer.
    There is a constant $C$ depending only on $k$ and $\epsilon$ such that $\unvd(D_k(n_1,n_2)) \leq (2^k+\epsilon)(n_1+n_2)+C$.
\end{theorem}

\begin{proof}
    For fixed constants $k,\epsilon$, simple computations show the existence of constants $C_0 \in \mathbb{N},\eta >0$ such that $\frac{\binom{C}{k}}{\binom{C(\frac{1}{2}-\eta)}{k}} \leq 2^k+\epsilon$ for every integer $C \geq C_0$.
    We let $C$ be the maximum between $C_0$ and $R_k\left(\left\lceil\frac{k}{\frac{1}{2}-\eta}\right\rceil,\left\lceil\frac{k}{\frac{1}{2}-\eta}\right\rceil\right)$.

    We first assume that $n_1,n_2 \geq \frac{C}{\eta(2^k+\epsilon)}$.
    Let $T$ be a tournament on $\lceil(2^k+\epsilon)n_1\rceil + \lceil (2^k+\epsilon)n_2\rceil+ C$ vertices.
    Consider a median order of $T$ and let $B_1$ be the first $\lceil(2^k+\epsilon) n_1\rceil$ vertices in this order, $B_2$ be the last $\lceil(2^k+\epsilon) n_2\rceil$ vertices,
    and $A$ be set of the remaining $C$ vertices in the middle.

    By~\ref{item:M2}, every vertex $v \in A$ is dominated by at least $\frac{1}{2}|B_1| - |A| \geq (\frac{1}{2}-\eta)|B_1|$ vertices in $B_1$ (using $|B_1| = \lceil(2^k+\epsilon)n_1\rceil \geq \frac{C}{\eta}=\frac{|A|}{\eta}$).
    Let $A' \subseteq A$ be any set of size $\left\lceil\frac{k}{\frac{1}{2}-\eta}\right\rceil$.
    By Lemma~\ref{lemma:kovaru_sos_turan_adapte}, there is a subset $A^\star$ of $A'$ of size $k$ such that $A^\star$ has at least $n_1$ common in-neighbours in $B_1$.
    Since this holds for every subset $A'$ of size $\left\lceil\frac{k}{\frac{1}{2}-\eta}\right\rceil$, and by choice of $C$, by Ramsey's Theorem there is a subset $\Tilde{A}$ of $A$ of size $\left\lceil\frac{k}{\frac{1}{2}-\eta}\right\rceil$ such that every subset of $k$ vertices in $\Tilde{A}$
    have at least $k$ common in-neighbours in $B_1$.
    By repeating the same reasoning, because every vertex in $A$ dominates at least $\frac{1}{2}|B_2| - |A| \geq (\frac{1}{2}-\eta)|B_2|$ vertices in $B_2$,
    we deduce that $\Tilde{A}$ contains a set $X$ of $k$ vertices with at least $n_2$ common out-neighbours in $B_2$.
    Since every $k$-subset of $\Tilde{A}$ has at least $n_1$ common in-neighbours in $B_1$, this gives the desired copy of $D_k(n_1,n_2)$.

    \medskip

    We now consider the case $n_1 < \frac{C}{\eta(2^k+\epsilon)}$ or $n_2 < \frac{C}{\eta(2^k+\epsilon)}$.
    We apply the previous case for $n_1,n_2$ replaced respectively by $n_1+\left\lceil \frac{C}{\eta(2^k+\epsilon)} \right\rceil$ and $n_2 + \left\lceil \frac{C}{\eta(2^k+\epsilon)} \right\rceil$.
    This shows that $\unvd(D_k(n_1,n_2)) \leq (2^k+\epsilon)\left(n_1+n_2+2\left(\frac{C}{\eta(2^k+\epsilon)}+1\right)\right)+C \leq (2^k+\epsilon)(n_1+n_2)+C'$ for $C'=C+2(2^k+\epsilon)\left(\frac{C}{\eta(2^k+\epsilon)}+1\right)$.
\end{proof}

While Proposition~\ref{prop:D0n} shows that this upper is tight when $n_1$ or $n_2$ is in $o(n)$, we do not know the precise asymptotic of $\unvd(D_k(n_1,n_2))$ when $n_1 \approx n_2$.

\subsection{\texorpdfstring{$1$}{1}-extensions of directed paths}

In this section, we give the exact value of $\unvd(D)$ when $D$ is a full $1$-extension of a directed path.

\begin{theorem}
Let $D$ be a full $1$-extension of a directed path with added vertex $s$.
If $s$ is a source or a sink, or if $D$ is $TT_3$, then $\unvd(D) = 2|V(D)|-2$. 
Otherwise, $\unvd(D) = 2|V(D)|-3$.
\end{theorem}
\begin{proof}
Observe that every orientation of $K_4$ contains $TT_3$, and that the directed triangle does not contain $TT_3$, so $\unvd(TT_3) = 4$. We now assume that $D$ is not $TT_3$.

Set $n=|V(D)|$, and let $s$ be a vertex of $D$ such that $D-s$ is a directed path $P$.
Let $V_1=N^-_D(s)$ and $V_2=N^+_D(s)$ and set $n_1=|V_1|$ and $n_2=|V_2|$.
Observe that $D = D_1(n_1,n_2)$. 
By Proposition~\ref{prop:unvd_star} we obtain $\unvd(D) \geq 2n-3$, and $\unvd(D) \geq 2n-2$ if $s$ is a sink or a source.

If $n_1=0$ or $n_2=0$, then $s$ is a source or a sink. Thus, by Proposition~\ref{prop:double}, we have $\unvd(D) \leq 2 \unvd(P) = 2n-2$.

It remains to show that $\unvd(D) \leq 2n-3$, assuming $n_1\geq 1$ and $n_2\geq 1$, and $\max \{n_1,n_2\} \geq 2$. By directional duality, assume without loss of generality that $n_1 \geq 2$.
Let $T$ be a tournament of order $2n-3$, and let $\sigma_0=(v_{-2n_1+1}, \dots , v_0, \dots ,  v_{2n_2-1})$ be a median order of $T$.
Set $A_1 = \{v_{-2n_1+1}, \dots , v_{-1}\}$ and $A_2 =  \{ v_0, \dots ,  v_{2n_2-1}\}$.
Let $B_1=N^-(v_{-1}) \cap A_1$ and $B_2=N^+(v_{-1}) \cap A_2$. Observe that, by \ref{item:M2}, $|B_2| \geq n_2$ and $|B_1| \geq n_1 - 1$. We first consider the case $|B_1| \geq n_1$.

\begin{claim}\label{claim:v0_dominated}
    If $|B_1| \geq n_1$, then $T$ contains a copy of $D$ with $s=v_{-1}$.
\end{claim}

\begin{proofclaim}
Let $I^-$ be the set of in-generators of $T\ind{B_1}$ and  $O^+$ be the set of vertices which are initial vertices of a directed path of order $n_2$ in $T\ind{B_2}$.
Note that $I^-\neq \emptyset$ and $O^+\neq \emptyset$. Moreover, by Lemma~\ref{lem:reduc-origin}, $(B_1\setminus I^-) \Ra I^-$ and $O^+ \Ra (B_2\setminus O^+)$.
Let $v_{i^-}$ be the vertex in $I^-$ with largest index and $v_{i^+}$ be the vertex in $O^+$ with smallest index.
 By Lemma~\ref{lem:origin}, $v_{i^-}$ is the terminus of a directed path $P^-$ of order $n_1$ in $T\ind{B_1}$ and  $v_{i^+}$ is the initial vertex of a directed path $P^+$ of order $n_2$ in $T\ind{B_2}$.
 Hence we may assume $v_{i^+}\ra v_{i^-}$, for otherwise the union of $v_{-1}$, $P^-$, $P^+$, the arcs between those paths and $v_{-1}$, and $v_{i^-}v_{i^+}$ is a copy of $D$ in $T$.

Assume first that $|B_2|= n_2$. Then, since $|A_2| = 2n_2$, $v_{-1}$ has exactly $n_2$ out-neighbour and $n_2$ in-neighbours in $A_2$.
Consequently $\sigma'=(v_{-2n_1+1}, \dots v_{-2}, v_0, \dots ,  v_{2n_2-1}, v_{-1})$ is also a median order of $T$.
By \ref{item:M4} for $\sigma'$, since $v_{i^+} \ra v_{i^-}$, there is a vertex $v_{\ell}$ in $\{v_{i^-}, \dots , v_{i^+}\}\setminus \{v_{-1}\}$ such that $v_{i^-} \ra v_{\ell} \ra v_{i^+}$.
Observe that $v_\ell \notin I^-$ by maximality of $i^-$ and  $v_\ell \notin (B_1\setminus I^-)$  because $(B_1\setminus I^-) \Ra \{v_{i^-}\}$.
Thus $v_{\ell}\notin B_1$, and similarly $v_{\ell}\notin B_2$. In particular, $v_\ell \notin (V(P^-) \cup V(P^+))$. If $v_{-1}\ra v_\ell$, remove the end-vertex from $P^+$, otherwise remove the initial vertex of $P^-$. In both cases, the resulting union of  $v_{-1}$, $P^-$, $P^+$ and the arcs between those paths and $v_{-1}$ is a copy of $D$ in $T$ (with $s= v_{-1}$) as desired.

Henceforth, we may assume that $|B_2| > n_2$. In particular it implies $|O^+|\geq 2$: to see this, take any Hamiltonian directed path of $T\ind{B_2}$ (which exists by Theorem~\ref{thm:redei}), then the two first vertices of this directed path are initial vertices of directed paths of order $n_2$ in $T\ind{B_2}$.
If there is a vertex $v_{\ell}$ in $\{v_{i^-}, \dots , v_{i^+}\}\setminus \{v_{-1}\}$ such that $v_{i^-} \ra v_{\ell} \ra v_{i^+}$, then as above there is a copy of $D$ in $T$ (with $s=v_{-1}$), so assume that there is no such $v_{\ell}$.
Hence, by \ref{item:M4} and \ref{item:M5}, $\Tilde{\sigma}= (v_{-2n_1+1}, \dots ,v_{i^--1},  v_{i^-+1}, \dots , v_{i^+},  v_{i^-},  v_{i^++1}, \ldots ,  v_{2n_2-1})$ is a median order of $T$.
Let $v_{j^+}$ be the vertex of smallest index in $O^+\setminus\{v_{i^+}\}$.
Let $Q^+$ be a directed path of order $n_2$ in $T\ind{B_2}$ with initial vertex $v_{j^+}$.
If $v_{i^-}\ra v_{j^+}$, then  the union of $v_{-1}$, $P^-$, $Q^+$, the arcs between those paths and $v_{-1}$, and $v_{i^-}v_{j^+}$ is a copy of $D$ in $T$. Henceforth we may assume $v_{j^+}\ra v_{i^-}$. 
By \ref{item:M4} for $\Tilde{\sigma}$, there is a vertex $v_p$ with $i^++1 \leq p \leq j^+-1$ such that $v_{i^-} \ra v_{p} \ra v_{j^+}$. 
Note that $v_p\neq v_{-1}$ because $i^+\geq 0$, and as above $v_p \notin B_1\cup B_2$.
As above, by removing either the end-vertex of $Q^+$ or the initial vertex of $P^-$, we get that the union of  $v_{-1}$, $P^-$, $Q^+$ and the arcs between those paths and $v_{-1}$ is a copy of $D$ in $T$.
\end{proofclaim}

Since Claim~\ref{claim:v0_dominated} settles the case $|B_1|\geq n_1$, and because $|B_1| \geq n_1-1$,  we now assume that $|B_1| = n_1 - 1$. 
Then by \ref{item:M4} and \ref{item:M5} applied on $\sigma_0$, we have that $\sigma_1 = (v_{-1},v_{-2n_1+1}, \dots, v_{-2} , v_0, \dots ,  v_{2n_2-1})$ is a median order of $T$. In particular this implies $v_{-2} \ra v_0$ and $v_{-1}\ra v_{-2n_1+1}$. Observe that $v_{-1},v_{-2n_1+1}, \dots, v_{-1}$ is then a Hamiltonian directed cycle of $T\ind{A_1}$.
Not also that if $|N^-(v_2) \cap A_1|\geq n_1$ we can apply Claim~\ref{claim:v0_dominated} with $\sigma_1,v_{-2}$ playing the roles of $\sigma_0,v_{-1}$ respectively.
By repeating this process, either we eventually find a copy of $D$ in $T$ with added vertex $s\in A_1$ or we ensure that $A_1 \Ra \{v_0\}$. Henceforth we assume that $A_1 \Ra \{v_0\}$.

The end of the proof is similar to the proof of Claim~\ref{claim:v0_dominated}, with $v_0$ playing the role of $v_{-1}$.
Let $B_1'= A_1$ and $B_2'=N^+(v_0) \cap A_2$. We have $|B_1'| = 2n_1 -1$ and, by \ref{item:M2} on $\sigma_0$, $|B_2'|\geq n_2$.
Since $v_{-1},v_{-2n_1+1}, \dots, v_{-1}$ is a directed cycle, in particular both $v_{-1}$ and $v_{-2}$ are the terminal vertices of a directed path of order $n_1$ in $B_1'$. 
Let ${O^+}'$ be the set of vertices which are initial vertex of a directed path of order $n_2$ in $T\ind{B_2'}$.
Note that ${O^+}'\neq \emptyset$ because $|B_2'|\geq n_2$ and, by Lemma~\ref{lem:reduc-origin}, ${O^+}' \Ra (B_2\setminus {O^+}')$.
Let $v_{k^+}$ be the vertex in ${O^+}'$ with smallest index.
Let $R^-$ be the directed path $v_{-n_1},\dots,v_{-1}$ and $R^+$ be a directed path of order $n_2$ in $T\ind{B_2}$ starting on $v_{k^+}$.

We may assume that $v_{k^+}\ra v_{-1}$, for otherwise $v_0$, $R^-$, $R^+$, the arcs between those paths and $v_0$, and $v_{-1}v_{k^+}$ is a copy of $D$ in $T$. 
If there is a vertex $v_{\ell}$ in $\{v_{1}, \dots , v_{k^+-1}\}$ such that $v_{-1} \ra v_{\ell} \ra v_{k^+}$, then by removing either the end-vertex of $R^+$ or the initial vertex of $R^-$, we get that the union of  $v_{0}$, $R^-$, $R^+$ and the arcs between those paths and $v_{0}$ is a copy of $D$ in $T$. 
This holds because $v_\ell \notin B_2'$ by choice of $v_{k^+}$ and $R^+$ is a directed path in $B_2'$.
Henceforth, we assume that there is no such $v_{\ell}$.
Therefore, by \ref{item:M4} and \ref{item:M5} applied on $\sigma_0$, $\hat{\sigma}= (v_{-2n_1+1}, \dots ,v_{-2},  v_{k^+},  v_{-1}, \dots , v_{k^+-1},   v_{k^++1}, \ldots ,  v_{2n_2-1})$ is a median order of $T$. In particular, $v_{-2} \ra v_{k^+}$. 
Let ${R^-}'$ be $v_{-n_1-1},\dots,v_{-2}$. 
We thus obtain that the union of $v_{0}$, ${R^-}'$, $R^+$ and the arcs between those paths and $v_{0}$ is a copy of $D$ in $T$.
\end{proof}

\bibliographystyle{abbrv}
\bibliography{refs}

\end{document}